\documentclass[12pt]{amsart}
\usepackage{graphicx}
\usepackage{tikz}
\usepackage{float}
\usepackage{caption}
\usepackage{amsmath}
\usepackage{amscd}
\usepackage{amsthm,enumerate}
\usepackage{amssymb}
\usepackage{setspace}
\usepackage{enumitem}
\usepackage{needspace}
  \usepackage{ytableau}

\usepackage[pagebackref=true]{hyperref}

\newtheorem{thm}{Theorem}[section]
\newtheorem{cor}[thm]{Corollary}
\newtheorem{lem}[thm]{Lemma}
\newtheorem{prop}[thm]{Proposition}
\theoremstyle{definition}
\newtheorem{defn}[thm]{Definition}
\newtheorem{exam}[thm]{Example}
\newtheorem{quest}[thm]{Question}
\newtheorem{rem}[thm]{Remark}

\numberwithin{equation}{section}
\DeclareMathOperator{\re}{Re}
\newcommand{\C}{\mathbb C}

\begin{document}

\title[Stability of well-covered Ferrers graphs]{Stability of independence polynomials of well-covered Ferrers graphs}
\author{Ta Quang Anh}
\address{Faculty of Mathematics and Informatics, Hanoi University of Science and Technology, 1 Dai Co Viet, Bach Mai, Hanoi, Vietnam.}
\email{quanganhta2710@gmail.com}
\author{Do Trong Hoang$^{*}$}
\address{Faculty of Mathematics and Informatics, Hanoi University of Science and Technology, 1 Dai Co Viet, Bach Mai, Hanoi, Vietnam.}
\email{hoang.dotrong@hust.edu.vn}

\let\thefootnote\relax\footnotetext{*Corresponding Author:  Do Trong Hoang -- {\tt hoang.dotrong@hust.edu.vn}}

\author{Nguyen Minh Nhat}
\address{Faculty of Mathematics and Informatics, Hanoi University of Science and Technology, 1 Dai Co Viet, Bach Mai, Hanoi, Vietnam.}
\email{nhatnguyeminh@gmail.com}


\date{}
\begin{abstract}
We prove that every root of the independence polynomial of a
well-covered Ferrers graph lies in the open unit disk centered at
$-1$, that is, $|x+1|<1$. Consequently, the polynomial is strictly
stable. We also exhibit a general Ferrers graph whose independence polynomial
is not stable.
\end{abstract}

\onehalfspacing
\maketitle

\section{Introduction}

Throughout this paper, all graphs are finite, simple, and undirected. For a graph   $G$, let $V(G)$ and $E(G)$ 
denote its vertex and edge sets, respectively.  
A subset of $V(G)$ consisting of pairwise non-adjacent vertices is  called an {\it independent set}. An independent set is \emph{maximal} if it is not
properly contained in another independent set.  It is \emph{maximum} if
it has the largest possible cardinality, denoted by $\alpha(G)$.
A graph is \emph{well-covered} if all its maximal independent sets have
the same cardinality.

 For each integer $k\geq 0$, let $s_k$
 denote the number of
independent sets of $G$ of cardinality $k$, with   $s_0=1$. The
\emph{independence polynomial} of $G$ is defined by 
\[
 I(G;x)=\sum_{k=0}^{\alpha(G)}s_k x^k.
\]

The roots of $I(G;x)$, known as \emph{independence roots},
have been extensively studied  (see, for example, \cite{BC, CS, LM17, WZ}). 
 Chudnovsky and Seymour~\cite{CS} prove that all independence roots   of     claw-free graph are real. In contrast,  independence  roots of   well-covered graphs need not be real, even for  trees (see  \cite{LM17}). 
Nevertheless, several results describe the distribution
of independence roots of well-covered graphs.     
Brown, Dilcher, and Nowakowski~\cite[Theorem~3.2]{BDN}
showed that if $G$ is a well-covered graph on $n\geq 1$ vertices,
then every independence root $x$ satisfies
$\frac{1}{n}\leq |x|\leq \alpha(G).$
On the other hand, Brown, Hickman, and Nowakowski~\cite{BHN}
proved that the independence roots of well-covered graphs
are dense in $\mathbb{C}$.

\begin{defn}
A nonzero   polynomial $f(x)$  is \emph{stable} if all its zeros lie
in the closed left half-plane $\{x\in\C:\re x\leq 0\}$. In other words, every root   $x\in \C$   of   $f(x)$  satisfies  $\re x\le 0$.
 The polynomial \(f(x)\) is called  \emph{strictly stable} if all its zeros lie in the open left half-plane, that is, if   $\re x<0$ for every root $x$ of $f(x)$. 
\end{defn}

A graph $G$ is called \emph{stable} if its independence polynomial
$I(G;x)$ is stable. Brown and Cameron~\cite{BC} proved that
every star $K_{1,n}$ is stable and posed the question of \emph{whether
all complete bipartite graphs are stable}.
Chen, Ning, and Tu~\cite{CNT} partially answered this question
by establishing the stability of $K_{2,n}$ and $K_{3,n}$
for every positive integer $n$.
They also showed that, for any fixed rational number $\ell>1$,
the graph $K_{m,\ell m}$ is nonstable whenever $m$ is
sufficiently large and $\ell m$ is an integer.
  
A related disk bound can be derived from the root correspondence
established by Levit and Mandrescu~\cite{LM08} for the corona
$G\circ K_1$, obtained by attaching a pendant vertex to each
vertex of $G$. Specifically, every root $x\neq -1$ of
$I(G\circ K_1;x)$ can be expressed as
$x=\frac{r}{1-r},$ 
where $r$ is a root of $I(G;x)$. If $G$ is stable, then
$\operatorname{Re}(r)\leq 0$ and $r\neq 0$, which implies
 $|x+1|=\frac{1}{|1-r|}<1.$ 
Since $-1$ also lies in this disk, every independence root
of $G\circ K_1$ lies in the open unit disk centered at $-1$.  In particular, the independence polynomial of   \(G\circ K_1\) is strictly stable.

  To investigate the locations of independence roots, it is
convenient to introduce the shifted polynomial
\[
\Phi_G(z)=I(G;z-1).
\]
Under the change of variables $z=x+1$, each root $x$ of
$I(G;x)$ corresponds to a root $z=x+1$ of $\Phi_G(z)$.
In particular, $G$ is stable if and only if every root $z$
of $\Phi_G(z)$ satisfies $\operatorname{Re}(z)\leq 1$.
Moreover, this transformation maps disks centered at $-1$
to disks centered at the origin, facilitating the root-location
analysis that follows.
 
Ferrers graphs form a natural class of bipartite graphs that
includes all complete bipartite graphs and has close connections
to commutative algebra. Corso and Nagel~\cite{CN} studied
their associated Ferrers ideals, proving that these ideals
have linear resolutions and computing the Hilbert series
of their quotient rings.
More precisely, let
$\lambda=(\lambda_1,\ldots,\lambda_n)$ be a partition with
$\lambda_1\geq\cdots\geq\lambda_n\geq 1$.
The Ferrers diagram of $\lambda$ consists of left-justified
rows of boxes of lengths $\lambda_1,\ldots,\lambda_n$.
The associated Ferrers graph $G_\lambda$ is the bipartite
graph whose adjacency structure is determined by these
row lengths, as described below.

\begin{defn}\label{def:ferrers} Let $\lambda=(\lambda_1,\ldots,\lambda_n)$ be a partition with
$m=\lambda_1\geq\cdots\geq\lambda_n\geq 1$. The \emph{Ferrers graph} associated with \(\lambda\)   is a bipartite graph $G_{\lambda}$  with vertex set $V(G_{\lambda}) = \{x_1,\ldots, x_n\}\cup \{y_1,\ldots, y_m\}$ and edge set 
$$E(G) = \{x_iy_j\mid 1\le i\le n, 1\le j\le \lambda_i\}.$$
Equivalently, the neighborhoods $N_{G_\lambda}(x_i)=\{y_1,\ldots,y_{\lambda_i}\}$, for $1\le i\le n$,  are nested. In particular, complete bipartite
graphs correspond to partitions with all parts equal.
\end{defn}

 The paper is organized as follows. In Section~2, we characterize
well-covered Ferrers graphs. In Section~3, we derive an explicit
formula for the independence polynomial of a Ferrers graph
from the Hilbert series of its associated quotient ring
(Lemma~\ref{lem:polynomial}). We then  
 prove  that every    independence root 
of  a well-covered
Ferrers graph lies
in the open disk $|x+1|<1$ (Theorem~\ref{thm:main}).
We also determine the sharp bound on the root modulus
and identify the rightmost root
(Proposition~\ref{prop:radius}).
Finally, we investigate Ferrers graphs obtained by adjoining
full columns to well-covered Ferrers graphs. We show that
adjoining one or two columns preserves strict stability,
whereas for any fixed $c\geq 3$, adjoining $c$ columns
produces nonstable graphs when the number of distinct
row lengths is sufficiently large
  (Theorem~\ref{prop:small-offset}).

\section{Well-covered Ferrers graphs}

The well-coveredness of Ferrers graphs has already been
characterized from a commutative algebraic perspective.
Let $G_\lambda$ be the Ferrers graph associated with a partition
$\lambda=(\lambda_1,\ldots,\lambda_n)$, where
$m=\lambda_1\geq\cdots\geq\lambda_n\geq 1$.
Its edge ideal in the polynomial ring
$R=k[x_1,\ldots,x_n,y_1,\ldots,y_m]$ is given by
\[
I_\lambda=(x_iy_j\mid 1\leq i\leq n,\ 1\leq j\leq\lambda_i).
\]
Since the complements of maximal independent sets are precisely
the minimal vertex covers, which correspond to the minimal
primes of $I_\lambda$, the graph $G_\lambda$ is well-covered
if and only if $I_\lambda$ is unmixed.
Corso and Nagel~\cite[Corollary~2.6]{CN} characterized
unmixed Ferrers ideals in terms of an anti-diagonal condition:
$n=m$, and every inside corner of the Ferrers diagram
lies on the line $p+q=m$.
This characterization was later extended to skew Ferrers
ideals by Hoang and Vu~\cite[Theorem~1.3]{HV}
using corner blocks. Motivated by these results, we provide an alternative characterization of well-covered Ferrers graphs.

  For the remainder of the paper, we adopt the following notation. 
Let $d_1>d_2>\cdots>d_t>0$ be the distinct parts of $\lambda$,
and let $r_j$ denote the multiplicity of $d_j$.   Thus, 
\[
\lambda=
(\underbrace{d_1,\ldots,d_1}_{r_1},
\underbrace{d_2,\ldots,d_2}_{r_2},
\ldots,
\underbrace{d_t,\ldots,d_t}_{r_t}),
\]
where $d_1=m$.  We refer to the group of $r_j$ rows of length
$d_j$ as the $j$-th \emph{row block} of $\lambda$ and denote by
$t(\lambda)=t$ the number of row blocks. 
Finally, set $p_0=0$ and  
 $p_j=\sum_{h=1}^j r_h$, for $1\leq j\leq t.$  In particular, $p_t=n$.  

For the associated Ferrers graph $G_\lambda$, let
$X_j=\{x_{p_{j-1}+1},\ldots,x_{p_j}\}$,
for  $1\leq j\leq t$, 
and $Y=\{y_1,\ldots,y_m\}$.
The following proposition characterizes the well-coveredness
of $G_\lambda$ in terms of the distinct row lengths $d_j$
and their multiplicities $r_j$.

\begin{prop}\label{prop:wc}
Let $G_\lambda$ be a Ferrers graph without isolated vertices. The following are equivalent:
\begin{enumerate}[label=\textup{(\roman*)}]
 \item $G_\lambda$ is well-covered.
 \item $ d_j=r_j+r_{j+1}+\cdots+r_t$, for all $j=1,\ldots,t$. 
  \end{enumerate}
\end{prop}

\begin{proof} 
For $j=2,\ldots,t$, by Corso and Nagel~\cite[Corollary~2.6]{CN},
the inside corner at the transition to row group $j$ has coordinates
$(p_{j-1},d_j)$. Thus the anti-diagonal condition is
$p_{j-1}+d_j=m$. Under $n=m$, this is equivalent to
\[
 d_j=n-p_{j-1}=\sum_{h=j}^t r_h,
 \qquad j=2,\ldots,t.
\]
For $j=1$, the same formula follows from $d_1=m=n=\sum_{h=1}^t r_h$.
Hence the unmixedness criterion is exactly condition~\textup{(ii)},
proving the equivalence.
\end{proof}

\begin{exam}\label{ex:4421}
Let $\lambda=(4,4,2,1)$ be the partition illustrated below.

\begin{center}
\begin{figure}[h!]
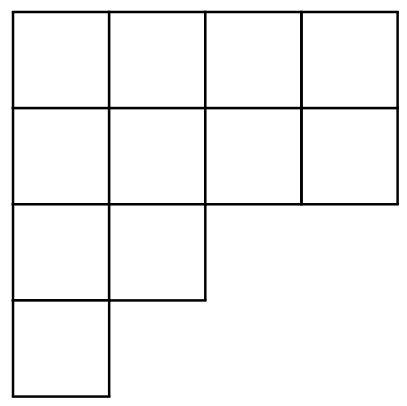

	     $\begin{ytableau}
 *(white) &*(white) &*(white) &*(white)    \\
 *(white) &*(white) &*(white) &*(white)   \\
  *(white) &*(white)  \\
   *(white)  
 \end{ytableau}$
 \caption{Well-covered Ferrers diagram}
\end{figure} 
\end{center} 
Here $m=n=4$,   and the distinct row lengths and their multiplicities
are $(d_1,d_2,d_3)=(4,2,1)$, and  $(r_1,r_2,r_3)=(2,1,1).$ 
Equivalently,
\[
 d_1=2+1+1=4,\qquad d_2=1+1=2,\qquad d_3=1.
\]
Proposition~\ref{prop:wc} therefore shows that $G_\lambda$ is
well-covered.   
\end{exam}

\section{Stability of well-covered Ferrers graphs}

In this section, we investigate the locations of independence
roots of well-covered Ferrers graphs and the stability of
Ferrers graphs obtained by adjoining full columns.
We begin by deriving an explicit formula for the independence
polynomial of an arbitrary Ferrers graph using the Hilbert
series of its associated quotient ring.

Recall that the quotient ring $R/I_\lambda$ inherits
the standard grading from $R$. For each integer $d\geq 0$,
let $(R/I_\lambda)_d$ denote its homogeneous component
of degree $d$. The Hilbert function of $R/I_\lambda$ is
defined by $H_{R/I_\lambda}(d)
=\dim_k (R/I_\lambda)_d,$  
and its Hilbert series is given by
\[
P_{R/I_\lambda}(q)
=\sum_{d\geq 0}H_{R/I_\lambda}(d)q^d.
\]

\begin{thm}[{\cite[Theorem~2.1]{CN}}] 
\label{thm:corso-nagel}  The Ferrers ideal $I_{\lambda}$ has   
a    2-linear minimal  free resolution. Then the Hilbert series of its quotient  $R/I_{\lambda}$ is  
\[   P_{R/I_\lambda}(q)
        =\frac{1}{(1-q)^m}
         +\frac{q}{(1-q)^{m+n+1}}
          \sum_{j=1}^{n}(1-q)^{\lambda_j+j}.
\]
\end{thm}

\begin{lem}\label{lem:polynomial}
   The independence polynomial of  Ferrers graph $G_\lambda$  is   
\begin{equation*} 
 I(G_\lambda;x)=(1+x)^m+
 x\sum_{i=1}^n(1+x)^{m+n-\lambda_i-i}.
\end{equation*}
Consequently, its shifted independence polynomial has the form
\begin{equation*} 
 \Phi_{G_\lambda}(z)=z^m+(z-1)
 \sum_{i=1}^n z^{m-i+n-\lambda_i}.
\end{equation*}
\end{lem}

\begin{proof}  
By \cite[Theorem 5.1.7]{BH}, the Hilbert series of the Ferrers ideal satisfies 
\[P_{R/I_\lambda}(q)
        =I\left(G_\lambda;\frac{q}{1-q}\right).
\]
Substituting  $q=\frac{x}{1+x}$, we obtain 
\[
\begin{aligned}
 I(G_\lambda;x)
 &=P_{R/I_\lambda}\!\left(\frac{x}{1+x}\right)\\
 &=(1+x)^m+
   \frac{x}{1+x}(1+x)^{m+n+1}
   \sum_{i=1}^n (1+x)^{-\lambda_i-i}\\
 &=(1+x)^m+
   x\sum_{i=1}^n(1+x)^{m-i+n-\lambda_i},
\end{aligned}
\]  
which establishes the first formula. The second follows immediately by substituting $x=z-1$.
\end{proof}

\begin{cor}\label{cor:wc-poly}
Let $G_\lambda$ be a well-covered Ferrers graph without isolated
vertices. Then  $m=n$, and   its shifted independence polynomial is
given by 
\begin{equation*}  
 \Phi_{G_\lambda}(z)=(t+1)z^n-\sum_{j=1}^t z^{n-r_j},
\end{equation*}
where  $r_1,\ldots,r_t$ are the sizes of the
row blocks of $\lambda$. 
\end{cor}

\begin{proof}
Since $G$ is well-covered, Proposition~\ref{prop:wc} gives $m=n$ and
\[
 d_j = r_j + r_{j+1} + \cdots + r_t = n - (r_1 + \cdots + r_{j-1})
\]
for every $j$. Note that  the sum $r_1 + \cdots + r_{j-1}$ is $0$ when $j=1$. By Lemma \ref{lem:polynomial}, 
\[
 \Phi_{G_{\lambda}}(z)   = z^n + (z-1)\sum_{i=1}^{n} z^{2n-i-\lambda_i}.
\]
   For $1\leq s\leq r_j$, the $s$-th row of the $j$-th row block
corresponds to row $i=p_{j-1}+s$ of the Ferrers diagram.  Every row in this group has length $\lambda_i = d_j $.  Thus, 
$i = (r_1 + \cdots + r_{j-1}) + s$ and 
$\lambda_i = n- (r_1 + \cdots + r_{j-1})$.  Hence, 
\[
 2n - i - \lambda_i
 = 2n - (r_1 + \cdots + r_{j-1}) - s - n + (r_1 + \cdots + r_{j-1})
 = n - s.
\]
Splitting the sum over the $t$ groups, we obtain
\begin{align*}
 \Phi(G;z)
 &= z^n + (z-1)\sum_{j=1}^{t}\sum_{s=1}^{r_j} z^{n-s}  = z^n + \sum_{j=1}^{t} \left(z^n - z^{n-r_j}\right) \\
 &= z^n + t z^n - \sum_{j=1}^{t} z^{n-r_j}  = (t+1)z^n - \sum_{j=1}^{t} z^{n-r_j}.
\end{align*}
\end{proof}

\begin{thm}[Rouch\'e's theorem]
\label{thm:rouche}
Let $\gamma$ denote a simple closed contour, and suppose that $F(z)$ and $H(z)$ are analytic inside and on $\gamma$. If $|H(z)|<|F(z)|$ at each point on $\gamma$ then $F(z)$ and $F(z)+H(z)$ have the same number of zeros, counting multiplicities, inside $\gamma$.
\end{thm}

\begin{thm}\label{thm:main}
Let $G_{\lambda}$ be a well-covered Ferrers graph without isolated vertices.
Every zero $x$ of $I(G_{\lambda};x)$ satisfies
   $|x+1|<1.$
  In particular, $I(G_{\lambda};x)$ is strictly stable, and hence stable.
\end{thm}

\begin{proof}
From  Corollary~\ref{cor:wc-poly}, let $\Phi_{G_{\lambda}}(z)=F(z)+H(z)$, where 
\[
 F(z)=(t+1)z^n,\qquad H(z)=-\sum_{j=1}^t z^{n-r_j}.
\]
  Since $1 \leq r_j \leq n$, both $F(z)$ and $H(z)$ are polynomials, and hence analytic on and inside  the unit circle $\gamma = \{z \in \mathbb{C} : |z| = 1\}$. On $\gamma$, we have:
\[
 |H(z)|\leq\sum_{j=1}^t|z|^{n-r_j}=t<t+1=|F(z)|.
\]
By Rouch\'e's theorem,   since $F$ has a zero of
multiplicity $n$ at the origin, $\Phi_{G_{\lambda}}(z)$ has exactly $n$ zeros in
$|z|<1$, counted with multiplicity. Since  $\Phi_{G_{\lambda}}(z)$ has
degree $n$,  all its zeros lie in $|z|<1$.

Therefore, if  $x$ is a zero of $I(G_\lambda;x)$, then $z=x+1$ is zero of $\Phi_{G_\lambda}$, so  
$\re x=\re z-1\leq |z|-1<0.$ 
Hence $I(G_\lambda;x)$ is strictly stable.
\end{proof}

\begin{rem}\label{rem:isolated}
The same conclusion holds if isolated vertices are allowed in the
definition of a Ferrers graph. Indeed, let $G^\circ$ be obtained by
deleting all $q$ isolated vertices. Every maximal independent set of
$G$ contains these isolated vertices, so $G$ is well-covered if and only if
$G^\circ$ is well-covered. Moreover,
\[
 I(G;x)=(1+x)^q I(G^\circ;x).
\]
The additional roots are all $-1$, which satisfy $|x+1|<1$.
If $G^\circ$ is empty, its independence polynomial is $1$, and the
same argument applies directly.
\end{rem}

Theorem~\ref{thm:main} gives the uniform bound $|x+1|<1$ for the
independence roots of every well-covered Ferrers graph without isolated vertices. We now sharpen
this bound for every graph in this class by determining the smallest closed disk
centered at $-1$ that contains all its roots, and show that the bound
is attained by a real root.

\begin{prop}\label{prop:radius}
Let $G_{\lambda}$ be a well-covered Ferrers graph without isolated vertices,
 whose row blocks have sizes  $r_1,\ldots,r_t$, where $r_j\geq1$ for all $1\leq j\leq t$.
Then there is a unique real number $p\in(0,1)$ satisfying
\begin{equation*}\label{eq:p}
    \sum_{j=1}^t p^{-r_j}=t+1.
\end{equation*}
Moreover,  every zero $x$ of $I(G_\lambda;x)$ satisfies
 $|x+1|\leq p$
 with 
 equality attained by a real zero. 
\end{prop}

\begin{proof}
Let define  the function
\[
    f(s)=\sum_{j=1}^t s^{-r_j}, \qquad s>0.
\]
Since  $f$ is continuous and strictly decreasing on
$(0,+\infty)$, with  
\[
    \lim_{s\to0^+}f(s)=+\infty
    \qquad\text{and}\qquad f(1)=t<t+1,
\]
  there exists a unique 
$p\in(0,1)$ such that $f(p)=t+1$.  

Now let $x$ be any zero of $I(G_\lambda;x)$
and set $z=x+1$.  By Corollary \ref{cor:wc-poly},  we have
\[
    \Phi_{G_{\lambda}}(z) =(t+1)z^n-\sum_{j=1}^t z^{n-r_j}.
\]
If $z=0$ then clearly  $|z|\le p$. Otherwise, dividing by $z^n$ gives   
  $\sum_{j=1}^t z^{-r_j}=t+1.$ 
  Taking absolute values and applying the
triangle inequality, we obtain  
\[
    f(p)=t+1
    =\left|\sum_{j=1}^t z^{-r_j}\right|
    \leq\sum_{j=1}^t |z|^{-r_j}
    =f(|z|).
\]
Since $f$ is strictly decreasing, it follows
that $|z|\leq p$,  or equivalently, $|x+1|\le p$. 
 
Finally, the defining equation for $p$ yields
\[
    \Phi_{G_\lambda}(p)
    =p^n\left(t+1-\sum_{j=1}^t p^{-r_j}\right)
    =0.
\]
Thus $x=p-1$ is a real zero of
$I(G_\lambda;x)$ satisfying $|x+1|=p$.
Hence the bound is sharp. 
\end{proof}

The well-covered assumption   in Theorem~\ref{thm:main} is essential. Indeed, the following  result 
constructs a family of Ferrers graphs with stable
independence polynomials. These graphs are obtained
from well-covered Ferrers graphs by adding a fixed
number of full columns, each corresponding to a
new vertex adjacent to every row vertex. 

\begin{thm}\label{prop:small-offset} Let $G_{\lambda}$ be a well-covered Ferrers graph. For each    positive integer  $c$, write $\lambda+c=(\lambda_1+c,\ldots, \lambda_n+c)$. 
Then the following statements hold.
\begin{enumerate}
    \item If $1\le c\le 2$, then      $G_{\lambda +c}$  
is stable. 
\item If $c\geq 3$ and $t(\lambda)\geq 3^c$, then $G_{\lambda+c}$
    is not stable.
\end{enumerate} 
\end{thm}

\begin{proof}   Since   $c>0$,   these two sets are maximal independent sets of different cardinalities, so  $G_{\lambda+c}$  is not well-covered. By Lemma \ref{lem:polynomial}, the  shifted independence polynomial of $G_{\lambda+c}$  is
\begin{eqnarray*}
    \Phi_{G_{\lambda+c}}(z) &=& z^{n+c}+ (z-1) \sum_{i=1}^{n}z^{2n+c-i-\lambda_i}.   
\end{eqnarray*}
 Let $t=t(\lambda)$,  
and  denote the distinct row lengths of $\lambda$ by
$d_1>\cdots>d_t$, with corresponding multiplicities
$r_1,\ldots,r_t$. Set 
$p_{j-1}=r_1+\cdots+r_{j-1}$ for $1\leq j\leq t$, where $p_0=0$.  Since $G_\lambda$ is
well-covered, Proposition~\ref{prop:wc} gives
$d_j=n-p_{j-1}$. In the extended diagram, a row $i=p_{j-1}+s$ in
 the  $j$-th row block, where $1\leq s\leq r_j$, has length $\lambda_i=d_j+c$. Thus the
exponent of this row in Lemma~\ref{lem:polynomial} is
\[
2n+c-i-(d_j+c)
 =2n-(p_{j-1}+s)-d_j
 =n-s.
\]

The added $c$ columns therefore cancel in the exponent. Summing the
terms from group $j$ and using
$(z-1)\sum_{s=1}^{r_j}z^{n-s}=z^n-z^{n-r_j}$, we get
\[
\begin{aligned}
\Phi_{G_{\lambda+c}}(z)
 &=z^{n+c}+(z-1)\sum_{j=1}^t\sum_{s=1}^{r_j}z^{n-s} =z^{n+c}+tz^n-\sum_{j=1}^t z^{n-r_j}\\
 &= z^n\left(z^c + t - \sum_{j=1}^t z^{-r_j}\right). 
\end{aligned}
\]

\vskip0.5em  
\noindent  (1): Suppose that $\Phi_{G_{\lambda + c}}(z)$ has a zero $z$ with
$\operatorname{Re}(z)\geq1$. Then $|z|\geq1$, and dividing
by $z^n$ gives
\[
    z^c+t=\sum_{j=1}^{t}z^{-r_j}.
\]
Since $r_j\geq1$ for every $j$, it follows that
\begin{equation*} 
    |z^c+t|
    \leq\sum_{j=1}^{t}|z|^{-r_j} \le \sum_{j=1}^{t}|z|^{-1} =\frac{t}{|z|}.
\end{equation*}
  
If $c=1$, then $|z+t|
    \geq\operatorname{Re}(z)+t
    \geq t+1
    >\frac{t}{|z|},$
a contradiction. 
Now suppose that $c=2$. Set  $a=\operatorname{Re}(z)$ and $s=|z|^2$.
Since  $a\geq1$ and  $s\geq1$, we have 
 $|z^2+t|^2=(s-t)^2+4ta^2.$  
It follows that
\[
\begin{aligned}
    |z|^2|z^2+t|^2-t^2
    &=s(s-t)^2+4tsa^2-t^2\\
    &\geq s(s-t)^2+4ts-t^2\\
    &=(s-1)(s-t)^2+s^2+2st >0.
\end{aligned}
\]
Thus $|z^2+t|>\frac{t}{|z|}$, again a contradiction. We conclude that every zero  $z$ of $\Phi_{G_{\lambda+c}}(z)$ satisfies
$\operatorname{Re}(z)<1$. Since
$\Phi_{G_{\lambda+c}}(z)=I(G_{\lambda+c};z-1)$, every zero $x$ of
$I(G_{\lambda+c};x)$ satisfies $\operatorname{Re}(x)<0$.
Hence $G_{\lambda + c}$ is stable for $1 \le c \le 2$.

\vskip0.5em  
\noindent  (2): Suppose that  $c\geq 3$ and  $t\ge 3^c$. Then $z$ is a
zero of $\Phi_{G_{\lambda+c}}(z)$ if and only if
\[
z^c+t-\sum_{j=1}^t z^{-r_j}=0.
\]
Set $z=t^{1/c}w$ and divide this equation by $t$. We obtain
\[
    w^c+1-\frac{1}{t}\sum_{j=1}^t(t^{1/c}w)^{-r_j}=0.
\]
Define  
\[
    F(w)=w^c+1,
    \qquad
    H_t(w)=-\frac{1}{t}
        \sum_{j=1}^{t}(t^{1/c}w)^{-r_j}.
\]
Then   $\Phi_{G_{\lambda+c}} (z)= tz^n  \left[F(w)+H_t(w)\right]$.   Observe that the zeros of  $F(w)$   are 
$w_k=e^{\pi(2k+1)i/c}$ for $0\le k\le c-1$, and are all simple. Moreover,  since $c\geq 3$, we have 
\[
    \operatorname{Re}w_0=\cos\frac{\pi}{c}
    \in\left[\frac12,1\right).
\]

For $k=1,\ldots,c-1$, we have
\[
|w_k-w_0|
 =\left|e^{(2k+1)\pi i/c}-e^{\pi i/c}\right|
 =\left|e^{2k\pi i/c}-1\right|
 =2\sin\frac{k\pi}{c}\geq 2\sin\frac{\pi}{c},
\] where equality holds if and only if $k=1$ or $k=c-1$.

Consider the open disk 
\[
    D=\{w\in\mathbb C:|w-w_0|<\frac{1}{2c}\}.
\]
We establish the following three claims:
\vskip0.5em \noindent 
{\it Claim 1.} $\overline{D}$ 
contains only one zero of $F(w)$, namely $w_0$. 
\vskip0.5em \noindent 
Indeed, recall that
$\sin x\geq 2x/\pi$ for all $x\in [0,\pi/2]$. Hence,  since $0<\pi/c\leq\pi/3<\pi/2$,   it follows that 
\[
    2\sin\frac{\pi}{c}
    \geq \frac4c
    > \frac1{2c}.
\]
Thus the distance from $w_0$ to every other zero of $F(w)$ exceeds
the radius of $\overline{D}$. Hence $w_0$ is the unique zero of $F(w)$ in
$\overline{D}$.

\vskip0.5em \noindent 
{\it Claim 2.}  ${\rm Re}(w)>\frac{1}{3}$ for all   $w\in D$. 
\vskip0.5em \noindent 
Indeed,  for every $w\in D$, we have $|w-w_0|<\frac{1}{2c}$.
Since $c\geq 3$, it follows that
\[
\begin{aligned}
    \operatorname{Re}(w)
    &\geq \operatorname{Re}(w_0)-|w-w_0|
    > \cos\frac{\pi}{c}-\frac1{2c}
     \geq \frac12-\frac16
     =\frac13.
\end{aligned}
\]
\vskip0.5em 

 By Claim 2,   $0\notin\overline{D}$. Since $F(w)$ is   polynomial   and $H_t(w)$  is analytic on $\mathbb C\backslash\{0\}$, both functions are analytic
on and inside $D$.

\vskip0.5em \noindent 
{\it Claim 3.}  $|H_t(w)| < |F(w)|$ on $\partial D$. 
\vskip0.5em \noindent  
Indeed, let $w\in\partial D$. By the reverse triangle inequality, we have  $|w| \geq |w_0|-|w-w_0| \geq 1- \frac{1}{2c}
\geq \frac{5}{6}$  on $\partial D$.  Since     $t\geq 3^c$ by the assumption, it follows that  
\[
    t^{1/c}|w|\geq 3\cdot\frac{5}{6}=\frac{5}{2}>1.
\]
Thus, using  $r_j\geq 1$ for every $j$, we obtain
\[
\begin{aligned}
    |H_t(w)|
    &\leq \frac1t\sum_{j=1}^t
       \frac{1}{(t^{1/c}|w|)^{r_j}} \leq \frac1t\sum_{j=1}^t\frac{1}{t^{1/c}|w|} = \frac{1}{t^{1/c}|w|}
      \leq \frac{2c}{3(2c-1)}.
\end{aligned}
\]

 On the other hand, since $w_0^c=-1$,   we have 
\[
\begin{aligned}
    F(w)
    &= w^c-w_0^c
     = \bigl(w_0+(w-w_0)\bigr)^c-w_0^c \\
    &= cw_0^{c-1}(w-w_0)
     + \sum_{k=2}^{c}\binom{c}{k}w_0^{c-k}(w-w_0)^k.
\end{aligned}
\]
Using the reverse triangle inequality, together with $|w_0|=1$ and
$|w-w_0|=\frac{1}{2c}$  on $\partial D$, we obtain
\[
\begin{aligned}
    |F(w)|
    &\geq \bigl|cw_0^{c-1}(w-w_0)\bigr|
      - \left|\sum_{k=2}^{c}
        \binom{c}{k}w_0^{c-k}(w-w_0)^k\right| \\
    &\geq \frac{1}{2}
      - \sum_{k=2}^{c}\binom{c}{k} \left(\frac{1}{2c}\right)^k. 
\end{aligned}
\]
 
For $2\leq k\leq c$, let
$b_k=\binom{c}{k}(\frac{1}{2c})^k$.  For $2 \le k < c$, $0<  c - k < c$ and $k + 1 \ge 3$,  
\[
    \frac{b_{k+1}}{b_k}
    = \frac{c-k}{k+1}\cdot \frac{1}{2c}
    = \frac{c-k}{c}\cdot\frac{1}{2(k+1)}
    < \frac{1}{6}.
\]
It follows that $b_{k+1}<b_k/6$  for all $2\le k< c$.   Iterating this inequality yields  $b_{k+1}< \frac{b_2}{6^{k-1}}$ for all $3\le k\le c$. Therefore,   
\[
    \sum_{k=2}^{c} b_k <   \sum_{k=2}^{c} \frac{b_2}{6^{k-1}} 
    = b_2\cdot  \frac{6}{5}\cdot \left(1- \frac{1}{6^{c-1}}\right) <b_2 \cdot \frac{6}{5} 
    = \binom{c}{2}\cdot \left(\frac{1}{2c}\right)^2 \cdot \frac{6}{5}
    = \frac{3(c-1)}{20c}.
\]
It follows that
\[
    |F(w)|
    > \frac12-\frac{3(c-1)}{20c}
    = \frac{7c+3}{20c} \ge \frac{2c}{3(2c-1)} \ge |H_t(w)|,
\] 
for all $w\in \partial D.$ 

By \emph{Claim 3} and  Rouch\'e's theorem, $F(w)+H_t(w)$ and $F(w)$ have the same number
of zeros in $D$, counted with multiplicity. Since $F(w)$ has   a unique zero $w_0$ in $D$, which  is simple, it follows that $F(w)+H_t(w)$ has exactly
one zero in $D$, denoted by $\xi \in D$.   
Moreover, $z=t^{1/c}\xi$ is a zero of $\Phi_{G_{\lambda+c}}(z)$.  By \emph{Claim 2},  $\operatorname{Re}(\xi)>1/3$. 
Hence  
\[
    \operatorname{Re}(z)
    = t^{1/c}\operatorname{Re}(\xi)
    > \frac{t^{1/c}}{3}
    \geq 1.
\]
Consequently, $G_{\lambda+c}$ is not stable, as required. 
\end{proof}

\section{Conclusion}

The preceding results show that well-coveredness is sufficient but not necessary, for the stability of independence polynomials of Ferrers graphs. Indeed, Theorem~\ref{thm:main} establishes that every well-covered Ferrers graph is stable, whereas
the examples above demonstrate the existence of both
stable Ferrers graphs that are not well-covered and
Ferrers graphs whose independence polynomials are not stable.

These observations naturally lead to the following question.

\begin{quest}\label{quest:general}
Let $\lambda=(\lambda_1,\ldots,\lambda_n)$ be a partition with
$m=\lambda_1\geq\cdots\geq\lambda_n\geq1$.
What are the necessary and sufficient conditions on $\lambda$
for the independence polynomial of $G_\lambda$ to be stable? 
\end{quest}

 \section*{Declarations}
 
 
 \noindent \textbf{Conflict of interest/Competing interests} \newline
 The authors declare that they have no competing interests\newline
 \noindent \textbf{Ethical approval and consent to participate}\newline
 \noindent Not applicable. \newline
 \noindent \textbf{Consent for publication}\newline
 \noindent Not applicable. \newline
 \noindent \textbf{Availability of data, code and materials}\newline
 \noindent Data sharing not applicable to this work as no data sets were
 generated or analyzed during the current study. \newline
 \noindent \textbf{Authors' contribution}\newline
 All authors have contributed equally to this work.\newline
 \noindent \textbf{Funding} \newline
 This research received no external funding. \newline


\begin{thebibliography}{99}

\bibitem{BC}
J.~I. Brown and B. Cameron,
\emph{On the stability of independence polynomials},
Electron. J. Combin. \textbf{25} (2018), no.~1, Paper No.~1.46.
 

  
\bibitem{BDN}
J.~I. Brown, K.~Dilcher, and R.~J. Nowakowski,
\emph{Roots of independence polynomials of well covered graphs}, J. Algebraic Combin.  \textbf{11} (2000), no.~3, 197--210.

\bibitem{BHN}
J.~I. Brown, C.~A. Hickman, and R.~J. Nowakowski, \emph{On the location of roots of independence polynomials},
 J. Algebraic Combin.  \textbf{19} (2004), no.~3, 273--282.



\bibitem{BH}
 W. Bruns and J. Herzog, {\it Cohen-Macaulay rings}, Cambridge Studies in Advanced Mathematics  39, Cambridge University Press, Cambridge, 1993.



\bibitem{CNT}
G. Chen, B. Ning, and J. Tu,
\emph{The stability of independence polynomials of complete bipartite graphs},
preprint (2025), arXiv:2505.24381. 

 	\bibitem{CS} M. Chudnovsky, P. Seymour, \emph{The roots of the independence
 		polynomial of a claw-free graph},  Journal of Combinatorial Theory,
 	Series B \textbf{97} (2007) 350--357.

 
\bibitem{CN}
A. Corso and U. Nagel,
\emph{Monomial and toric ideals associated to Ferrers graphs},
Trans. Amer. Math. Soc. \textbf{361} (2009), no.~3, 1371--1395. 



\bibitem{HV} D. T. Hoang and T. Vu, 
{\em Unmixed and sequentially Cohen-Macaulay property of skew tableau ideals},   J. Algebr. Comb.  {\bf 62} (2025), 26.


\bibitem{LM08}
V.~E. Levit and E. Mandrescu,
On the roots of independence polynomials of almost all very well-covered graphs,
\emph{Discrete Appl. Math.} \textbf{156} (2008), no.~4, 478--491.


\bibitem{LM17} V.~E.~Levit, E.~Mandrescu, \emph{On clique corona graphs--a
 		short survey}, Matem\'{a}tica Contempor\^{a}nea \textbf{45} (2017) 68--76.

\bibitem{WZ} Y.~Wang, B.-X.~Zhu, \emph{On the unimodality of independence polynomials of some
graphs}, European J. Combin. {\bf 32} (2011) 10--20. 
 

\end{thebibliography}
\end{document}